\documentclass[12pt]{article} \usepackage{verbatim}
\usepackage{graphicx}
\usepackage{hyperref,url}
\usepackage{amsthm}
\newcommand{\nocopyright}{
No Copyright\thanks{
The authors hereby waive all copyright
and related or neighboring rights to this work,
and dedicate it to the public domain.
This applies worldwide.
}}
\title{Division by four}
\author{Peter G. Doyle \and Cecil Qiu}
\date{Version dated 6 April 2015
\\ \nocopyright
}

\newtheorem{theorem}{Theorem}

\newtheorem{prop}[theorem]{Proposition}

\newcommand{\proofstart}{\noindent {\bf Proof.\ }}
\newcommand{\proofend}{$\quad \qed$}
\newcommand{\mathproofend}{\quad \qed}

\newcommand{\union}{\cup}
\newcommand{\cross}{\times}

\newcommand{\implies}{\;\Longrightarrow\;}

\newcommand{\seteq}{\asymp}
\newcommand{\setleq}{\preceq}

\newcommand{\setll}{\ll}

\newcommand{\goesto}{\rightarrow}

\newcommand{\kw}{\mathbf}

\newcommand{\good}{\kw{good}}
\newcommand{\bad}{\kw{bad}}

\begin{document}

\maketitle

\begin{abstract}
Write $A \setleq B$ if there is an injection from $A$ to $B$,
and $A \seteq B$ if there is a bijection.
We give a simple proof
that for finite $n$,
$n \cross A \setleq n \cross B$ implies $A \setleq B$.
From the Cantor-Bernstein theorem it then follows
that $n \cross A \seteq n \cross B$ implies $A \seteq B$.
These results have a long and tangled history,
of which this paper is meant to be the culmination.

\end{abstract}
\centerline{\emph{For John}}

\section{The gist of it}

{\bf To show:}
If there is a one-to-one map from $4 \cross A$ to $4 \cross B$
(which need not hit all of the range $4 \cross B$),
then there is a one-to-one map from $A$ to $B$.

\noindent
{\bf A word to the wise:} Check out what Rich Schwartz has to say in
\cite{schwartz:four}.

\proofstart
Think of $4 \cross B$ as a deck of cards where for each $x$ in $B$
there are cards of rank $x$ in each of the four suits
spades, hearts, gems, clubs.
Note that while we use the word `rank', in this game all ranks will
be worth the same:
Who is to say that a king is worth more or less than a queen?

Think of $A$ as a set of racks, where each rack has four spots
to hold cards,
and think of $4 \cross A$ as the set of all the spots in all the racks.

Think of a one-to-one map from $4 \cross A$ to $4 \cross B$
as a way to fill the racks with cards,
so that all the spots have cards,
though some of the cards may not have been dealt out and are still in the deck.

Name the four spots in each rack by the four suits as they come in bridge,
with spades on the left, then hearts, gems, clubs.
Call a spade `good' if it is in the spades spot of its rack,
and `bad' if not.

Do these two rounds in turn.
(As you read what is to come, look on to where we have worked out
a case in point.)

{\bf Shape Up}:
If a rack has at least one bad spade and no good spade,
take the spade that is most to the left and swap it to the spades
spot so that it is now good.
Do these swaps all at the same time in all of the hands,
so as not to have to choose which swap to do first.

{\bf Ship Out}:
Each bad spade has a good spade in its rack,
thanks to the Shape Up round.
Swap each bad spade for the card whose rank is that of the good spade
in its rack, and whose suit is that of its spot.
To see how this might go, say that in some rack the queen of spades
is in the spades spot, while the jack of spades is in the hearts spot.
In this case we should swap the jack of spades for the queen of hearts.
(Take care not to swap it for the jack of hearts!)
Note that some spades may need to be swapped
with cards that were left in the deck,
but this is fine.
Do all these swaps at once, for all the bad spades in all the racks.
This works since no two bad spades want to swap with the same card.

{\bf Note.}  If you want, you can make it
so that when there is more than one bad spade in a rack,
you ship out just the one that is most to the left.

Now shape up, ship out, shape up, ship out, and so on.
At the end of time there will be no bad spades to be seen.
(Not that all bad spades will have shaped up, or been put back in the deck:
Some may be shipped from spot to spot through all of time.)
Not all the cards in the spades spots need be spades,
but no card to the right of the spades spot is a spade.
So if we pay no heed to the spades spots,
we see cards of three suits set out in racks with three spots each,
which shows a one-to-one map from $3 \cross A$ to $3 \cross B$.

You see how this goes, right?
We do a new pass, and get rid of the hearts from the last two spots in
each rack.
Then one last pass and the clubs spots have just clubs in them.
So the clubs show us a one-to-one map from the set $A$ of racks to the
set $B$ of ranks.
Done.

Is this clear?
It's true that we have
left some things for you to think through on your own.
You might want to look at
\cite{schwartz:four},
where Rich Schwartz has put in things that we have left out.

\newpage

Here's the first pass in a case where all the cards have been dealt out.
Note that in this case we could stop right here and use the spades to
match $A$ with $B$, but that will not work when $A$ and $B$ get big.

\footnotesize
\begin{verbatim}
Start:
 4g    6h    Qg    8g    9h   *Qs*   4c    Ag    6c   *4s*
 Jh    Ah    9c    8h   *As*   Tc    Tg    5h    Qc   *Js*
 Kc   *6s*   4h    6g   *Ts*  *9s*   Jc    Kg   *8s*   8c 
 5c    5g   *Ks*  *5s*   Th    Jg    Ac    Qh    9g    Kh 

Shape up:
 4g   *6s*  *Ks*  *5s*  *As*  *Qs*   4c    Ag   *8s*  *4s*
 Jh    Ah    9c    8h    9h    Tc    Tg    5h    Qc   *Js*
 Kc    6h    4h    6g   *Ts*  *9s*   Jc    Kg    6c    8c 
 5c    5g    Qg    8g    Th    Jg    Ac    Qh    9g    Kh 

Ship out:
 4g   *6s*  *Ks*  *5s*  *As*  *Qs*   4c   *Ts*  *8s*  *4s*
 Jh    Ah    9c    8h    9h    Tc    Tg    5h    Qc   |4h|
 Kc    6h   *Js*   6g   |Ag|  |Qg|   Jc    Kg    6c    8c 
 5c    5g   *9s*   8g    Th    Jg    Ac    Qh    9g    Kh 

Ship out:
 4g   *6s*  *Ks*  *5s*  *As*  *Qs*   4c   *Ts*  *8s*  *4s*
 Jh    Ah    9c    8h    9h    Tc    Tg    5h    Qc   |4h|
*9s*   6h   |Kg|   6g   |Ag|  |Qg|   Jc   *Js*   6c    8c 
 5c    5g   |Kc|   8g    Th    Jg    Ac    Qh    9g    Kh 

Shape up:
*9s*  *6s*  *Ks*  *5s*  *As*  *Qs*   4c   *Ts*  *8s*  *4s*
 Jh    Ah    9c    8h    9h    Tc    Tg    5h    Qc   |4h|
 4g    6h   |Kg|   6g   |Ag|  |Qg|   Jc   *Js*   6c    8c 
 5c    5g   |Kc|   8g    Th    Jg    Ac    Qh    9g    Kh 

Ship out:
*9s*  *6s*  *Ks*  *5s*  *As*  *Qs*   4c   *Ts*  *8s*  *4s*
 Jh    Ah    9c    8h    9h    Tc   *Js*   5h    Qc   |4h|
 4g    6h   |Kg|   6g   |Ag|  |Qg|   Jc   |Tg|   6c    8c 
 5c    5g   |Kc|   8g    Th    Jg    Ac    Qh    9g    Kh 

Shape up:
*9s*  *6s*  *Ks*  *5s*  *As*  *Qs*  *Js*  *Ts*  *8s*  *4s*
 Jh    Ah    9c    8h    9h    Tc    4c    5h    Qc   |4h|
 4g    6h   |Kg|   6g   |Ag|  |Qg|   Jc   |Tg|   6c    8c 
 5c    5g   |Kc|   8g    Th    Jg    Ac    Qh    9g    Kh 

\end{verbatim}
\normalsize

\section{Discussion and plan} \label{discuss}

\subsection{Division by any finite number}

Write $A \setleq B$ (`$A$ is less than or equal to $B$')
if there is an injection from $A$ to $B$.
Write $A \seteq B$
(`$A$ equals $B$', with apologies to the equals sign)
if there is a bijection.

The method we've described for dividing by four
works fine for any finite $n$,
so we have:

\begin{theorem} \label{thleq}
For any finite $n$,
$n \cross A \setleq n \cross B$ implies $A \setleq B$.
\end{theorem}

From the Cantor-Bernstein theorem
(Prop. \ref{cb} below)
we then get
\begin{theorem} \label{theq}
For any finite $n$,
$n \cross A \seteq n \cross B$ implies $A \seteq B$.
\end{theorem}

As an application of the Theorem \ref{thleq},
Lindenbaum and Tarski
(cf. 
\cite[p. 305]{lindenbaumTarski:ensembles},
\cite[Theorem 13]{tarski:cancellation})
proved the following:
\begin{theorem} \label{euclid}
If
$m \cross A \seteq n \cross B$ where $\gcd(m,n)=1$,
then there is some $R$
so that $A \seteq n \cross R$ and $B \seteq m \cross R$.
\end{theorem}

We reproduce Tarski's proof in Section \ref{sec:division} below.

Combining Theorems \ref{theq} and \ref{euclid}
yields this omnibus result for division of an equality
\cite[Corollary 14]{tarski:cancellation}:

\begin{theorem}
If
$m \cross A \seteq n \cross B$ where $\gcd(m,n)=d$,
then there is some $R$
so that $A \seteq (n/d) \cross R$ and $B \seteq (m/d) \cross R$.
\end{theorem}

\subsection{Pan Galactic Division}
We call the shape-up-or-ship-out algorithm for eliminating
bad spades \emph{Shipshaping}.
As we've seen, Shipshaping is the basis for a division algorithm
that we'll call
\emph{Pan Galactic Long Division}.
As the name suggests, there is another algorithm called
\emph{Pan Galactic Short Division},
which we'll come to presently.
`Pan Galactic' indicates that we think this is the `right way'
to approach division,
and some definite fraction of intelligent life forms in the universe
will have discovered it.
Though if this really is the right way to divide, 
there should be no need for this Pan Galactic puffery,
we should just call these algorithms \emph{Long Division} and
\emph{Short Division}.
Which is what we will do.

\subsection{What is needed for the proof?}
Shipshaping and its associated division procedures are effective (well-defined,
explicit, canonical, equivariant, \ldots),
and do not require the well-ordering principle or any other form
of the axiom of choice.
Nor do we need the axiom of power set,
which is perhaps even more suspect than the axiom of choice.
In fact we don't even need the axiom of infinity,
in essence because if there are no infinite sets that
all difficulties vanish.
Still weaker systems would suffice:
It would be great to know
just how strong a theory is needed.

\subsection{Whack it in half, twice}
Division by 2 is easy (cf. Section \ref{two}),
and hence so is division by 4:
\[
4 \cross A \setleq 4 \cross B \implies 2 \cross A \setleq 2 \cross B
\implies A \setleq B
.
\]
We made it hard for ourselves in order to show a method
that works for all $n$.
It would have been more natural to take $n=3$,
which is the crucial test case for division.
If you can divide by 2, and hence by 4,
there is no guarantee that you can divide by 3;
whereas if you can divide by 3, you can divide by any $n$.
This is not meant as a formal statement, it's what you might call a Thesis.
We chose $n=4$ instead of $n=3$
because there are four suits in a standard deck of cards,
and because there is already a paper called `Division by three'
\cite{conwaydoyle:three},
which this paper is meant to supersede.

\subsection{Plan} \label{plan}

In Section \ref{two} we discuss division by two, and explain why
it is fundamentally simpler than the general case.
In Section \ref{short} we introduce Short Division.
In Section \ref{pgs} we take a short break to play Pan Galactic Solitaire.
In Section \ref{cancel} we reproduce classical results
on subtraction and division,
so that this work can stand on its own as the definitive resource
for these results,
and as preparation for Section \ref{timing},
where we discuss how long
the Long and Short Division algorithms take to run.
In Section \ref{history}
we discuss the tangled history of division.
In Section \ref{vale}
we wrap up with some platitudes.

\section{Division by two} \label{two}
Why is division by two easy?
The flippant answer is that $2-1=1$.
Take a look at Conway and Doyle
\cite{conwaydoyle:three}
to see one manifestation of this.
Here is how this shows up in the context of Shipshaping.

The reason Shipshaping has a Shape Up round is that for $n>2$,
there can be more than one bad spade.
When $n=2$ there can't be more than one bad spade.
In light of this, we can leave out the Shape Up rounds.
When a bad spade lands in the hearts spot of a rack
with a heart in the spades spot, we just leave it there.
It's probably easier to understand this if we give up the good-bad
distinction, and just say that the rule is that when both cards in a rack
are spades, we ship out the spade in the hearts spot.
At the end of time, there will be at most one spade in each rack,
so in the Long Division setup there will be at least one heart;
assigning to each rack the rank of the rightmost heart gives
us an injection from racks to ranks.

In this approach to division by 2,
we find that there is no need to worry about doing everything in
lockstep.
We can do the Ship Out steps in any order we please,
without organizing the action into rounds.
As long as any two-spade rack eventually gets attended to,
we always arrive at the same final configuration of cards.
This is the kind of consideration that typically plays an important
role in the discussion of distributed computation,
where you want the result not to depend on accidents of what happens first.
It's quite the opposite of what concerns us here, where,
in order to keep everything canonical, we can't simply say,
`Do these steps in any old order.'
Without the axiom of choice, everything has to be done
synchronously.

Now in fact the original Shipshaping algorithm
works fine as an asynchronous algorithm for any $n$,
but the limiting configuation will depend on the order in which swaps are
carried out,
so the result won't be canonical.
More to the point, without the Axiom of Choice
we can't just say,
`Do these operations in whatever order you please.'
In contrast to the real world, where the ability to do everything
synchronously would be a blessing,
for us it is an absolute necessity.

\section{Short Division} \label{short}

In Short Division we reverse the roles of $A$ and $B$, so that
$A$ is the set of ranks and $B$ the set of racks.
An injection from $4 \cross A$ to $4 \cross B$ now shows a way to deal
out all the cards, leaving no cards in the deck, though some of
the spots in the racks may remain empty.
We do just one round of Shipshaping.
This works just as before, the only new twist being that if the
spades spot is empty when we come to shape up a bad spade,
we simply move the spade
over into the spades spot.
Since now all the cards have been dealt out, we don't ever
have occasion to swap with a card still in the deck.

When $A$ is finite, all the spades will shape up,
and show a bijection from ranks ($A$) to racks ($B$).

When $A$ is infinite,
some of the bad spades may get `lost',
having been shipped out again and again
without ever shaping up.
These lost spades will each have passed through
an infinite sequence of spots,
and all these infinite sequences will be disjoint.
We use these sequences to hide the lost spades,
as follows.

Let $A_\good$ denote the ranks of the good spades at the end of the game,
and $B_\good \seteq A_\good$ the racks where they have landed.
$A_\bad = A - A_\good$ is the the set of ranks of the lost spades.
To each element of $A_\bad$ there is an infinite chain of spots in
$3 \cross B_\good$, and these chains are all disjoint.
Using these disjoint chains we can use the usual `Hilbert Hotel' scheme
to define a bijection between
$A_\bad \union 3 \cross B_\good$ and $3 \cross B_\good$,
i.e. we make $3 \cross B_\good$ `swallow' $A_\bad$.
But if $3 \cross B_\good$ can swallow $A_\bad$ then so can $B_\good$
(cf. Proposition \ref{swallow} below):
$A_\bad \union B_\good \seteq B_\good$.
So
\[
A  = A_\bad  \union A_\good \seteq A_\bad \union B_\good \seteq B_\good
\setleq B
.
\]

{\bf Note.}  To make apparent the sequences of cards accumulated by the lost
spades,
we can modify the game
by making stacks of cards accumulate under the spades,
to record their progress.
The Shape Up round is unchanged,
except that we swap the whole spade stack into the
spades spot.
In the Ship Out round, when a spade ships out it takes its stack with it,
and places it on top of the card it was meant to swap with.
Spades that eventually shape up will have finite piles
beneath them,
but lost spades will accumulate
an infinite stack of cards.

\section{Pan Galactic Solitaire} \label{pgs}

Let us take a short break to play Pan Galactic Solitaire.

Deal the standard 52-card deck out in four rows of 13 cards.
The columns represent the racks, with spots in each rack 
labelled spades, hearts, diamonds, clubs going from top to bottom.
(Though as you'll see it makes no difference how we label the spots,
this game is suit-symmetric.)
The object is to `fix' each column
so that the ranks of the four cards are equal
and each card is in the suit-appropriate spot.
We move one card at a time.
There is no shaping up;
shipping out swaps are allowed based on any suit, not just spades.
So for example if in some column the 3 of hearts is in the hearts spot
and the 6 of hearts is in the diamonds spot,
we may swap the 6 of hearts for the 3 of diamonds.

We don't know a good strategy for this game.
Computer simulations show that various simple strategies yield a probability
of winning of about 1.2 per cent.
Younger players find this frustratingly low, and either play with a smaller
deck or allow various kinds of cheating.
Older players are not put off by the challenge,
and at least one has played the game often enough to have won twice.
Though because the game recovers robustly from an occasional error,
he cannot be certain that he won these games fair and square.

A couple of apps have been written to implement this game on the computer.
In some versions if you click on the 6 of hearts (as in the
example above)
the app locates the 3 of diamonds for you and carries out the swap.
This makes the game go much faster,
but it is much less fun to play than if you must locate and
click on the 3 of diamonds,
or better yet, drag the 6 of hearts over onto the 3 of diamonds.
One theory as to why the automated version of the game is less fun to play
is that the faster you can play the game,
the more frequently you lose.

This game is called Pan Galactic Solitaire from the conviction that something
like it will have occurred to a definite fraction of all life forms
that have discovered Short and Long Division.

\section{Cancellation laws} \label{cancel}

In this section we reproduce basic results about cancellation.
These results are all cribbed from
Tarski
\cite{tarski:cancellation},
though the notation is new and (we hope) improved.

To simplify notation, 
write $+$ for disjoint union and $n A$ for $n \cross A$.
(This abbreviation is long overdue.)
$A-B$ will mean set theoretic complement, with the implication
that $B$ is a subset of $A$.

The results below guarantee the existence of certain injections and bijections,
and the proofs are backed by algorithms.
For finite sets these can all be made to run snappily.

\subsection{Subtraction laws}

The only ingredients here are the Cantor-Bernstein construction
and the closely related Hilbert Hotel construction.

We start with Cantor-Bernstein.

\begin{prop}[Cantor-Bernstein] \label{cb}
\[
A \setleq B
\;\land\;
B \setleq A
\implies
A \seteq B
.
\]
\end{prop}

\proofstart
For a proof,
draw the picture and follow your nose
(cf.\ \cite{conwaydoyle:three}).

Here is a version of the proof emphasizing that the desired
bijection is the pointwise limit of a sequence of finitely-computable
bijections, which will be important to us in
Section \ref{timing} below.

It suffices to show that from an injection
\[
f: A \goesto B \subset A
\]
we can get a bijection.

Say we have any function
$f:A \goesto B$,
not necessarily injective
(this relaxation of the conditions is useful,
so that we can think of finite examples).
Every $a \in A$  makes a Valentine card.
To start, every $a \in A-B$ gives their Valentine to $f(a) \in B$.
Any $b \in B$ that gets a Valentine then gives their own Valentine to
$f(b)$.
Repeat this procedure ad infinitum,
and at the end of the day
every $b \in B$
has at least one Valentine,
and every Valentine has a well-defined location
(it has moved at most once).
Let $g$ associate to $a \in A$ the $b \in B$ that has $a$'s
Valentine.
$g$ is always a surjection,
and if $f$ is an injection, $g$ is a bijection.

Here's a variation on the proof.
Again each $a \in A$ makes a Valentine,
only this time every $a \in A$ gives their Valentine to $f(a)$.
Now any $b \in B$ that has no Valentine demands its Valentine back.
Repeat this clawing-back ad infinitum,
and at the end of the day,
every Valentine has a well-defined location,
and if $f$ was injective, or more generally if $f$ was finite-to-one,
every $b \in B$ has a Valentine,
and we're done as before.
The twist here is that if $f$ is not finite-to-one,
at the end of the day some $b$'s may be left without a Valentine.
So they demand their Valentine back, continuing a transfinite
chain of clawings-back.
We may or may not be comfortable concluding that after some transfinite
time, every $b \in B$ will have a Valentine.
For present purposes we needn't worry about this, since
when $f$ is injective the fuss is all over after at most $\omega$ steps.
\proofend

\noindent
{\bf Notation.}
Write $A \setll B$
(`$A$ is swallowed by $B$' or `$A$ hides in $B$')
if $A + B \setleq B$.
By Cantor-Bernstein this is equivalent to $A + B \seteq B$.
Another very useful equivalent condition is that there exist
disjoint infinite chains inside $B$,
one for each element of $A$.

By repeated swallowing we have:

\begin{prop} \label{multiswallow}
For any $n$, if
\[
A_i \setll B
,\;
i=0,\ldots,n-1
\]
then
\[
A_0 + \ldots + A_{n-1} \setll B
.
\mathproofend
\]
\end{prop}

Here are two close relatives of Cantor-Bernstein, proved by the same
back-and-forth construction.

\begin{prop}
\[
A + C \setleq B + C
\implies
A - A_0 \setleq B
,
\]
where $A_0 \setll C$.
\proofend
\end{prop}

\begin{prop}
\[
A + C \seteq B + C
\implies
A - A_0 \seteq B - B_0
,
\]
where $A_0,B_0 \setll C$.
\proofend
\end{prop}

\begin{prop}
\[
A + C \setleq B + 2C
\implies A \setleq B + C
\]
\end{prop}

\proofstart
\[
A+C \setleq B + C + C
\implies A - A_0 \setleq B + C
\]
with $A_0 \setll C$.
So
\[
A = (A - A_0) + A_0 \setleq B + C + A_0 \setleq B + C
.
\mathproofend
\]

\begin{prop} \label{mnineq}
For finite $m < n$,
\[
A + m C \setleq B + n C
\implies
A \setleq B + (n-m) C
.
\]
\end{prop}

\proofstart
The proof is by induction, and we can get the number of recursive steps
down to $O(\log(m))$.
\proofend

From this we get

\begin{prop} \label{preswallow}
For $n \geq 1$
\[
A + n C \setleq B + n C
\implies
A + C \setleq B + C
.
\mathproofend
\]
\end{prop}

From Cantor-Bernstein we then get
\begin{prop}
For $n \geq 1$
\[
A + n C \seteq B + n C
\implies
A + C \seteq B + C
.
\mathproofend
\]
\end{prop}

Here's the key result for Short Division:

\begin{prop} \label{swallow}
For $n \geq 1$
\[
A \setll n C
\implies
A \setll C
.
\]
\end{prop}

\proofstart
This is the special case of Proposition \ref{preswallow} when $B$ is empty.
Here we redo the proof in this special case,
in preparation for the timing considerations of Section \ref{timing}.

Think of $nC$ as a deck of cards, where $(i,c)$ represents a card of suit
$i$ and rank $c$.
Since $A \setll nC$,
there are disjoint infinite chains
\[
s_a:\omega \goesto nC
,\;a \in A
.
\]
Let $\alpha(a)$ to be the smallest
$i$ such that $s_a$ contains infinitely many cards of suit $i$,
and let $\rho_{a}$ be the sequence of ranks of those cards.
(For future reference, note that we could trim this down to the
ranks of those cards of suit $\alpha(a)$ that come after the
last card of a lower-numbered suit.)

Let
\[
A_i = \{a \in A:\alpha(a)=i\}
.
\]
The infinite chains
\[
\rho_{a}: \omega \goesto C
,\;
a \in A_i
\]
are disjoint, so
\[
A_i \setll C
\]
So by Proposition \ref{multiswallow},
\[
A = A_0 + \ldots + A_{n-1} \setll C
.
\]
The required injection from $A+C$ to $C$ is obtained as a composition
of injections $f_i$ which map $A+C$ onto $A-A_i + C$, leaving $A-A_i$ fixed.
\proofend

Here's a result
that will be handy when we come to the
Euclidean algorithm.

\begin{prop} \label{handy}
For $m<n$
\[
A + m C \seteq n C
\implies
A + E \seteq (n-m) C
,
\]
where $E \setll nC$, and hence $E \setll C$.
\end{prop}

\proofstart
From above we have
\[
A \setleq (n-m) C
.
\]
Write
\[
A + E \seteq (n-m)C
,
\]
so that
\[
A + E + m C \seteq n C
.
\]
Since $A + mC \seteq n C$ this gives
\[
E + nC \seteq nC
\implies
E \setll nC
\implies
E \setll C
.
\mathproofend
\]

Finally,
here's the result Conway and Doyle needed to make their division method
work.

\begin{prop}
For $n \geq 1$,
\[
n A \setleq n B
\;\land\;
B \setleq A
\implies
A \setleq B
,
\]
and hence by Cantor-Bernstein $A \seteq B$.
\end{prop}

\proofstart
Write
\[
A \seteq B + C
.
\]
From
\[
nA \setleq nB
\]
and
\[
A = B+C
\]
we get
\[
nB+nC \setleq nB
,\]
i.e.
\[
nC \setll nB
.
\]
But
\[
nC \setll n B \implies C \setll n B \implies C \setll B
,
\]
so
\[
A \seteq B + C \setleq B
.
\mathproofend
\]

\subsection{Division laws} \label{sec:division}

Finally we come to division.
We've already proved Theorems \ref{thleq} and \ref{theq}.
All that remains now is Theorem \ref{euclid}.

{\bf Proof of Theorem \ref{euclid}.}
As you would expect, 
the proof is a manifestation of
the Euclidean algorithm.
If $m=1$ we are done.
Otherwise we will use the usual recursion.

\[
mB \setleq n B \seteq m A
.
\]
Using division we get
\[
B \setleq A
.
\]
Write
\[
A \seteq B + C
.
\]
\[
mB+mC \seteq mA \seteq n B
.
\]
From Lemma \ref{handy} we have
\[
mC + E \seteq (n-m) B
,
\]
with $E \setll B$, and hence $E \setll A$.
We think of $E$ as `practically empty'.
If it were actually empty, we'd recur using $C$ in place of $A$.
Ditto if we knew that $E \setll C$.
As it is we use $C+E$ in the recursion,
hoping that this amounts to practically the same thing.
\[
m(C+E) \seteq mC + E + (m-1)E
\seteq (n-m)B + (m-1)E
\seteq (n-m)B
.
\]
So by induction, for some $R$ we have
\[
C+E \seteq (n-m)R
,
\]
\[
B \seteq m R
.
\]
Since $E \setll A$,
\[
A \seteq A + E
\seteq B + C + E
\seteq
m R + (n-m) R
\seteq n R
\]
Done.

We won't undertake to determine
the best possible running time here.
But in order to make sure it requires at most a logarithmic number
of divisions,
we will want to check that for the recursion
we can subtract any multiple $km$ from $n$
as along as $km < n$.
Here's the argument again, in abbreviated form.
\[
km B \setleq nB \seteq m A
;
\]
\[
kB \setleq A
;
\]
\[
A \seteq kB+C
;
\]
\[
mkB+mC \seteq mA \seteq n B
;
\]
\[
mC+E \seteq (n-km)B
,\;
E \setll B \setleq A
;
\]
\[
m(C+E) \seteq mC+E+(m-1)E \seteq (n-km)B + (m-1)E \seteq (n-km)B
;
\]
\[
C+E \seteq (n-km)R
,\;
B \seteq mR
;
\]
\[
A \seteq A+E \seteq kB+C+E \seteq kmR + (n-km)R \seteq n R
.
\mathproofend
\]

\section{How long does it take?} \label{timing}

\subsection{The finite case}

In Shipshaping, every swap puts at least one card in a spot where it
will stay for good, so the number of swaps is at most $n|A|$:
This holds both in the Long Division setup where $nA$ is the set of spots,
and in
the Short Division setup where $nA$ is the set of cards.
For a distributed process where all Shape Up and Ship Out rounds take
place simultaneously
the number of rounds could still be nearly this big, despite the
parallelism,
because one bad spade could snake its way through nearly all the non-spades
spots.
If we simulate this distributed process on a machine with one processor,
the running time will be $O(n |A|)$
(or maybe a little longer,
depending on how much you charge for various operations).

{\bf Note.} While $|B|$ might be as large as $n|A|$, 
nothing requires us to allocate storage for all $n|B|$ cards (in Long Division)
or spots (in Short Divison).

For finite $A$,
in Short Division
all spades will shape up after one pass of Shipshaping,
and show an injection from ranks to racks.
So the running time for Short Division is $O(n|A|)$,
running either as a distributed process or on a single processor.
This is as good as we could hope for.

Still for finite $A$,
Long Division with the naive recursion takes $n-1$ passes of Shipshaping.
If we divide by 2 whenever an intermediate value of $n$ is even,
we can get this down to at most $O(\log(n))$ passes.
The number of suits remaining 
gets halved at least once every other pass,
so the total number of swaps over all rounds of Shipshaping will be at most
\[
|A|(n + n + n/2 + n/2 + n/4 + n/4 + \ldots)
=
4n|A|
.
\]
Hence the total time for Long Division is $O(n|A|)$,
running either as a
distributed process or on a single-processor machine.
This is the same order as for Short Division, though Short Division
will win out when you look at explicit bounds.

\subsection{The infinite case}

For $A$ infinite,
to talk about running times we will need a notion of transfinite
synchronous distributed computation.
The general idea is to support taking the kind of pointwise limits
that show up in Shipshaping,
where in the limit the contents of each spot is well defined,
as is the goodness (but not the location) of each spade.
(See \ref{timing} below for more specifics.)

One round of Shipshaping will take time $\omega$.
For Long Division sped up as in the finite case so that as to take
$O(\log(n))$ passes,
the running time will be $\omega \cdot O(\log(n))$.

{\bf Note.}
Here and throughout,
we'll be rounding running times down to the nearest limit ordinal,
so as not to have to worry about some finite number of post-processing
steps.

For real speed, we'll want to use Short Division.
The swallowing step can be implemented by a recursion which,
like Long Division, can be sped up to take $O(\log(n))$ passes.
Though the number of passes is on the same order as for Long Division,
this can still be considered an improvement,
to the extent that Shipshaping is more complicated than swallowing.

The big advantage of Short Division stems from the fact that
the swallowing can be sped up to run in time $\omega$,
in essence by running the steps of the recursion in tandem.
And the swallowing can be configured to
run simultaneously with the Shipshaping.
Sped up in this way, Short Division can be done in time at most $\omega$.
We discuss this further in Section \ref{timing} below.

By contrast,
we've never found a way to run the division stages of Long Division
in tandem.

\subsection{Dividing a bijection}
To divide a bijection to get a bijection,
we can simultaneously compute injections
each way, and then combine them using the Cantor-Bernstein
construction in additional time $\omega$.
(Cf. Proposition \ref{cb}.)
So if dividing an injection takes time $\omega$,
dividing a bijection takes time at most $\omega \cdot 2$.

Can this be improved upon?
It is tempting to start running 
the Cantor-Bernstein algorithm using partial information about the two
injections being computed, but we haven't made this work.
The case to look at first is $n=2$,
where Long Division is all done after one pass of Shipshaping.

\subsection{Speeding up Long Division}

A division algorithm takes as input an injection $f_n:nA \goesto nB$,
and produces as output an injection $f_1:A \goesto B$.
In the Long Division setup, where all the spots are filled though
not all the cards need have been dealt out, one pass
of Shipshaping gets us from $f_n$ to $f_{n-1}$,
or more generally,
from $f_{n}$ to $f_{n(k-1)/k}$, for any $k$ dividing $n$.
In particular, when $n$ is even one pass gets us $f_{n/2}$.
This allows us to get 
from $f_n$ to $f_1$ in $O(\log(n))$ passes.
As one Shipshaping pass takes time at most $\omega$, this makes for
a total running time of at most $\omega \cdot O(\log(n))$.

Various tricks can be used to cut down the number of passes
in Long Division.
We can run Shipshaping using any divisor of $n$ we please,
or better yet, run it simultaneously for all divisors.
Knowing $f_m$ for as many values of $m$ as possible
can be useful because if we know
injections from $m_k A$ to $m_k B$ then we can paste them together to get
an injection from $m A$ to $m B$ for any positive linear combination
$m = \sum_k r_k m_k$.
This pasting takes only finite time,
which in this context counts as no time at all.
Combining these observations we can shave down the number of Shipshaping
passes needed,
and in the process we observe some intriguing phenomena.
But we can never get the number of passes
below $\lceil \log_2 n \rceil$, which is
at most
a factor of two better than
what we achieve with the naive method of dividing by 2 when any intermediate
$n$ is even.

\subsection{Speeding up Short Division}

For real speed, we use Short Division.
In this setup, all the cards are dealt
out, though not all the spots need be filled.
As the Shipshaping algorithm runs,
we observe a steadily decreasing collection of bad spades,
together with steadily lengthening disjoint sequences in $(n-\{0\}) \cross B$
telling the spots through which these bad spades have passed.
In the limit, the bad spades that remain have wandered forever,
and they index disjoint injective sequences in $(n-\{0\}) \cross B$.

The proof of Proposition \ref{preswallow} offers a recursive algorithm for
hiding the lost bad spades in with the good, which when sped up in the by-now
usual way requires $O(\log(n))$ recursive `passes', which if carried out
one after the other, give us a total running time of
$\omega \cdot O(\log(n))$.
(We've silently absorbed the $\omega$ coming from the
single pass of Shipshaping.)
This is on the same order as what we achieved with Long Division.

Having to do only one round of Shipshaping could be viewed as an improvement,
on the grounds that a Shipshaping is more complicated than a swallowing pass.
But as we've stated before, the real advantage of Short Division 
come from the fact that we can run all the swallowing passes in
tandem with each other, and with the Shipshaping algorithm.

Here's roughly how it works.
As we run the Shipshaping pass for Short Division,
the set of bad spades decreases, while the sequences of spots they have
visited steadily increases.
As these preliminary results trickle in,
we can be computing
how we propose
to hide the shrinking set of bad spades among the growing set of good spades.
Eventually every bad spade knows where it will go,
as does every good spade.
In the limit we have an injection from spades, and hence ranks, to racks.
The whole procedure is done after only $\omega$ steps
(though as always we reserve the right to do
some fixed finite number of post-processing steps).

The main thing missing here is how we end up hiding the lost spades,
and how we compute this.
The fundamental idea is to trim the sequence of spots visited by
a lost spade down to the subsequence
consisting of all the hearts visited,
then all the diamonds after the last heart,
then all the spades after the last diamond.
These trimmed sequences can be computed in tandem with Shipshaping,
so they are ready for use after time $\omega$.
(See Section \ref{chips} below.)

Now because these sequences are increasing
we can determine the limiting suit of all these sequences
by running along them, keeping track of the suit,
which eventually stops changing.
After another $\omega$ steps we've computed what in the notation
of Proposition \ref{swallow} was the function $\alpha$.
Then we can distinguish the lost spades according to this limiting value.
We first hide those that limit with hearts, then diamonds, then clubs.
This hiding takes only $O(n)$ post-processing steps,
which we disregard, for a total running time of $\omega \cdot 2$:
$\omega$ to compute and trim the tracks of the lost spades;
$\omega$ to determine the limiting suits of the lost spades.

To get the running time down to $\omega$ is trickier.
The idea is that we start using the trimmed sequences before we are done
computing them.
The injection we compute will be different from that just described,
because of artifacts associated to the 
cards in the trimmed sequences that come before those of the limiting suit
(e.g. a finite number of hearts coming before an infinite number of diamonds).
We omit the details.

\subsection{Chipshaping} \label{chips}

Here is a variation on Shipshaping that incorporates 
the simultaneous determination of the trimmed sequences of spots that the
lost spades have passed through.

In this variation, we begin by placing a poker chip on top of each spade,
which will serve as its representative during the Shape Up and Ship Out rounds.
(If the racks are tilted, we'll have to lean the chip against the card,
maybe it would be better to think of the cards laid out in rows
as in Pan Galactic Solitaire.)
When we shape up a chip we move the card and chip together,
but when we ship out we move only the chip.
We add a third round called Trim:  If a card has a chip on it, and the card
is a spot to the left of the spot for its suit (e.g. a club in a hearts spot),
we leave the chip where it is, and swap the card underneath the chip
to where it belongs, meaning the spot for its suit in
the rack having the spade of its rank in the spades spot.
(If the chip has just moved, this is where the chip came from;
in any case it is someplace the chip has visited before.)
We repeat the Trim round until no more moves are possible, meaning that no
card that is topped by a chip is `above its station'.
This takes only a finite number of rounds (at most equal to the number of
rounds Ship Out rounds we have done).
Then we continue:  Shape Up, Ship Out, Trim, Shape Up, Ship Out, Trim, Trim, Shape Up, Ship Out, Trim, Trim, Trim,\ldots

If you try this, you'll see that it is really quite nice, though it is
annoying that we might have to wait
through an ever larger number of Trim rounds.
An alternative is to do just one Trim round, but this entails---well, try it,
and you'll see.

In any event, we aren't quite computing everything we'll need in order to
see an injection from $A$ to $B$ in the limit after only $\omega$ steps.
For that, it seems that
we might have to go beyond what you can conveniently do with just
the original deck of cards and some chips.
Like, say, add some local memory for pointers,
maybe in the form of stickers affixed to cards or racks.
There's nothing wrong with this from the point of view of distributed
synchronous computation.
It just won't be as much fun.

\subsection{Transfinite synchronous distributed computation} \label{trans}

We're using here a loose notion of transfinite computation.
Roughly speaking, we're imagining that we have a processor for each element of
$A$ and $B$, or maybe (for convenience) for each element of $n \cross A$
and $n \cross B$.
Each processor has some finite set of flags,
a finite set of registers that can store the name
of an element of $A$ or $B$, and a finite set of ports.
The processors can communicate by sending messages to a designated
port of another processor;
if two processors simultaneously attempt
to send to the same port of another processor,
the whole computation crashes.
We allow processors to set the flag of another processor;
if two or more processors simultaneously
try set the same flag, there is no conflict.
Certain flags and registers may be designated as suitable for reading at limit
times $\omega, \omega \cdot 2, \ldots$.
If the state of one of these flags or the contents of one of these registers
fails to take on a limiting value, the whole computation crashes.

This is as precise as we want to get here.
Really, we're hoping to find that a suitable formulation has already
been made and explored.
It would be great if there is essentially one suitable formulation,
but we are by no means certain of this.

\section{Some history} \label{history}

Here's a brief rundown on the history of Theorems \ref{thleq} and \ref{theq}.

Theorems \ref{thleq} and \ref{theq}
follow easily from the well-ordering principle,
since then $n \cross A \seteq A$ when $A$ is infinite.
The well-ordering principle is a consequence of the power set axiom
and the axiom of choice (which without the power set axiom may take
several inequivalent forms).
The effective proofs we are interested in don't use either axiom.

Briefly put, 
Bernstein stated Theorem 2 in 1905 and gave a proof for $n=2$,
but nobody could make sense of what he had written about extending the
proof to the general case.
In the 20's Lindenbaum and Tarski found proofs of Theorems 2 and 1
but didn't publish them.
Lindenbaum died and Tarski forgot how their proofs went.
In the 40's Tarski found and published two new proofs of Theorem 1.
In the 90's Conway and Doyle found a proof,
after finally peeking at Tarski
\cite{tarski:cancellation};
based on what Tarski had written, they decided that their proof was
probably essentially the same as Lindenbaum and Tarski's lost proof.

For the gory details,
here is Tarski 
\cite{tarski:cancellation}, from 1949.
(We'll change his notation slightly to match that used here.)

\begin{quotation}
Theorem 2 for $n=2$ was first proved by F. Bernstein
\cite[pp.\ 122 ff.]{bernstein:mengenlehre};
for the general case Bernstein gave only a rough outline of a proof,
the understanding of which presents some difficulties.
Another very elegant proof of Theorem 2 in the case $n=2$
was published later by W. Sierpinski
\cite{sierpinski:two};
and a proof in the general case was found, but not published,
by the late A. Lindenbaum
\cite[p. 305]{lindenbaumTarski:ensembles}.
Theorem 1 --- from which Theorem 2 can obviously be derived by means of the
Cantor-Bernstein equivalence theorem --- was first obtained for $n=2$
by myself,
and then extended to the general case by Lindenbaum
\cite[p. 305]{lindenbaumTarski:ensembles};
the proofs, however, were not published.
Recently Sierpinski
\cite{sierpinski:leq}
has published a proof of Theorem 1 for $n=2$.

A few years ago I found two different proofs of Theorem 1
(and hence also, indirectly, of Theorem 2).
\ldots\ 
The second proof is just the one which I should like to present in this paper.
It is in a sense an extension of the original proof given by me for $n=2$,
and is undoubtedly related to Lindenbaum's proof for the general case.
Unfortunately, I am not in a position to state how close this relation is.
The only facts concerning Lindenbaum's proof which I clearly remember
are the following:
the proof was based on a weaker though related result previously obtained by me,
which will be given below \ldots;
the idea used in an important part of the proof was rather similar to the one
used by Sierpinski in the above-mentioned proof of Theorem 2 for $n=2$.
Both of these facts apply as well to the proof I am going to outline.
On the other hand, my proof will be based upon a lemma \ldots, which seems
to be a new result and which may present some interest in itself;
it is, however, by no means excluded that the proof of this
lemma could have been easily obtained by analyzing Lindenbaum's argument
\end{quotation}

Observe Tarski's delicate description of Bernstein's attempts at the
general case.
Conway and Doyle
\cite{conwaydoyle:three}
were less circumspect:
\begin{quote}
We are not aware of anyone other than Bernstein himself
who ever claimed to understand the argument.
\end{quote}

Around 2010 Arie Hinkis claimed to understand Bernstein's argument
(see \cite{hinkis}).
Peter proposed to Cecil that he look into this claim for his
Darmouth senior thesis project.
Between 2011 to 2013 we spent many long hours
trying to understand Hinkis's retelling of
Bernstein's proof,
and exploring variations on it.
In the end, we hit upon the proofs given here.
These proofs
are very much in the spirit of Bernstein's approach,
but not close enough that we believe
that Bernstein knew anything very like this,
and we remain skeptical that Bernstein ever knew a correct proof.
There are very many possible variations of Bernstein's general idea,
a lot of which seem to \emph{almost} work.
Still, these two proofs can be seen as a vindication of Bernstein's approach.
Our finding them owes everything to Hinkis's faith in Bernstein.

\section{Valediction} \label{vale}

Tarski
\cite[p.\ 78]{tarski:cancellation}
wrote that `an effective proof of [Theorems \ref{thleq} and \ref{theq}]
is by no means simple.'
We hope that you will disagree,
if not from this treatment then from Rich Schwartz's
enchanting presentation in
\cite{schwartz:four}.

The reason that division can wind up seeming simple to us is that
combinatorial arguments and algorithms are second nature to us now.
A different way of looking at the problem makes its apparent difficulties
disappear.

For an even more persuasive example of this, consider the Cantor-Bernstein
theorem,
and the confusion that accompanied it when it was new.
From a modern perspective,
you just combine the digraphs associated to the two injections
and follow your nose.

After a century of combinatorics and computer science,
it's easy to understand how Pan Galactic Division works.
Can we hope that someday folks will marvel at how hard it was to
discover it?

\bibliography{four}
\bibliographystyle{hplain}
\end{document}